\documentstyle[12pt]{article}

 \input{amssym}

   \def\bZ{{\Bbb Z}}
   \def\bN{{\Bbb N}}

   \def\bF{{\Bbb F}}
   \def\bk{{l\!k}}

\def\hH{{H^{\frak H}_\ast }} 
\def\hC{{C^{\frak H}_\ast}} 

\def\oC{ {\cal C}} 
\def\oA{ {\cal A}} 
\def\oB{ {\cal B}} 

\def\oE{ {\cal BE}} 
\def\oO{ {\cal O}} 
\def\End{ {\cal E}nd} 

\def\W{ \mbox{\rm W}} 

\def\Hom{ \mbox{\rm Hom}} 

\def\C{\mbox{\bf  C}} 

\def\D{{\bf \Delta }}  

\def\GM{\bk  \mbox{-\rm \bf GM}} 

 \def\DGM{\bk   \mbox{-\rm \bf DGM}} 

 \def\Alg{\mbox{-\rm \bf ALG}} 

\def\Tot{\mbox{\rm Tot}} 

\begin{document}

\vspace{3cm }

\centerline{ \large \bf  Operadic Hochschild chain complex and  free loop
spaces.}

\vspace{3mm}

\centerline{ \small by}
\vspace{3mm}

\centerline{ {\bf {\sc David Chataur  and Jean-Claude Thomas}}}

\vspace{ 0,3 cm}

\vspace{ 0,2 cm}

      \noindent {\bf Abstract.} {\small We construct, for any algebra $A$ over
an operad $\oO$, an  Hochschild chain
complex,
$\hC(\oO, A)$ which is also an $\oO$-algebra. This Hochschild chain complex
coincides with the usual one, whenever $A$ is a
commutative differential graded algebra. Let $X$ is  a simply connected space,
$N^\ast (-)$ be the singular cochain functor, $
X^{S^1}$ be the free loop space,
$\oC_\infty$ be a cofibrant replacement of commutative operad and $M_X$ a
$\oC_\infty $-cofibrant model of $X$. We prove that The operadic
chain complex
$\hC(\oC_\infty, M_X)$ is quasi-isomorphic to  $N^\ast (X^{S^1})$ as a
$\oC_\infty$-algebra. In particular, for any prime
field of coefficients this identifies the  action of the large  Steenrod algebra
on the Hochschild homomology
$\hH \left (N^\ast (X)\right)$ with the usual Steenrod operations on   $H^\ast
\left(X^{S^1}\right)$.}

\vspace{ 0,2 cm}
\rule{16 cm}{0,2mm}

\vspace{ 0,3 cm}


\vspace{ 0,3 cm}


       \vspace{5mm}

      {\bf Introduction.}

As illustrated by the fundamental result of Mandell, \cite{[Ma1]},
$E_\infty$-algebras are the good tools for the study of
the homotopy theory of topological spaces. Indeed, for a prime field of coefficients  the homology of a
$E_\infty$-algebra is an unstable algebra over the large Steenrod algebra. This
last property originates with the previous works of Dold, \cite{[Do]}, May
\cite{[Ma]} and has been recently extensively studied by the first
author  and Livernet
\cite{[C-L]}.

In this paper we develop an application of the homotopy theory of
$E_\infty$-algebras to the study  of the free loop space,
$X^{S^1}$ on a space $X$. For this purpose, we  consider the
almost free operad, denoted $\oC_\infty$, which is a cofibrant
model of the commutative operad. A quotient of this operad, also
denoted $\oC_\infty$,  has been studied by Kadeisvili,
\cite{[Ka]}, Kontsevich, \cite{[Ko]} and Getzler-Jones,
\cite{[G-J]}. It appears that the category of
$\oC_\infty$-algebras is a closed model category which is very
similar to the category of commutative and associative
differential graded algebras.

\vspace{3mm}

It follows from Theorem C and Proposition 4.2-2:

{\bf Main Theorem:}{\it We denote by $N^\ast(-)$ the normalized
singular cochain functor with coefficients in an arbitrary
commutative ring $\bk$ and by $H^\ast(-)= H(N^\ast(-))$ the
functor of singular cohomology. Let $X$ be a 1-connected space, if
each $H^i(X)$ is finitely generated then there exist natural
equivalences of $\oC_\infty$-algebras between $N^*(X^{S^1})$ and
the operadic Hochschild complex of the $\oC_\infty$-algebra
$N^\ast(X)$.}

\vspace{3mm}

 An  associative $\oC_\infty$-algebra is a particular case of  strongly homotopy algebra (see
3.3). Our result is, in this special  case, a  substancial
improvement of the results proved in \cite{[BT1]} and \cite{[BT2]}
since the structure of unstable algebra over the Steenrod algebra
but also all secondary cohomological operations are preserved
\cite{[C-L]}.

\vspace{3mm}

All the paper is devoted to the proof of this result. The required
knowledge about operads  is presented in section 1.  In section 2
we define the operad $\oC_\infty$ in relation with the
Barratt-Eccles operad, extensively studied by Berger and Fresse,
\cite{[B-F]}. In section 3, we construct for any operad $\oO$ in
the category of differential graded modules the operadic
Hochschild complex of an $\oO$-algebra, $A$. When $A$ is supposed
to be an associative and commutative differential graded algebra
we compare the  operadic Hochschild complex of an  $\oO$-algebra
with the usual Hochschild complex. We also study the particular
case when $A$ is an  almost free $\oO$-algebra (Theorem B). In the
last section we end the proof of our main result (Theorem C).

\vspace{5mm }

{\bf 1.  Backgrounds about algebras over an operad.}

\vspace{3mm }

\noindent{\bf 1.1 Notation.} We denote by $\GM$ (resp. $\DGM$) the
category of graded modules (resp. of differential graded modules).
We also consider the forgetfull functor:
$$
\# : \DGM \to \GM \,, \quad (V,d) \mapsto (V,d)_\#=V\,.
$$
We are mainly concerned by the following
categories :

 $\D $,  the simplicial category of finite ordered  sets with objects
$[n]=\{0,1,...n\} $ and non
decreassing maps

$  \C ^{\D}$, the  category  of cosimplicial objects and cosimplicial maps of
$\C$: $\underline X = \left( \{\underline
X^n\}_{n\geq 0}  , d ^i, s^i\right) $

$  \C ^{\D ^{op}}$, the  category  of simplicial objects and simplicial maps of
$\C$:  $\underline X = \left( \{\underline
X_n\}_{n\geq 0}  , d _i, s_i\right) $

\vspace{3mm }

\noindent{\bf 1.2 Operads. } Recall from \cite{[G-K]},
\cite{[G-J]} and \cite{[K-M]}, that an operad $\oO$ is defined in
any symmetric monoidal category $\C$ as a sequence of left $\bk
[\Sigma_i]$-modules (where $\Sigma_i$ is the symmetric group)
$\oO(i)\,, i \geq 0$, with {\it composition products}
$$
\oO(n)\otimes \oO(i_1)\otimes \oO(i_2)\otimes ... \otimes \oO(i_n) \to
\oO(i_1+i_2+...+i_n) \,, \quad x_0 \otimes x_1\otimes
...\otimes x_n \mapsto x_0(x_1,x_2,...,x_n)
$$
which are equivariant,  associative and with a unit. Homomorphisms
of operads are defined in an obvious way.

The category of operads is a closed model category
\cite{[B-M]}, \cite{[H]}, where the weak equivalences are
the quasi-isomorphisms  and the fibrations are the surjections.

The universal example of operad is the {\it endomorphim
operad}. Let $(V,d_V)$
be a differential graded module,
$\End _V$ is an operad in $\DGM$ such that:
$$
\End _V (n) = \Hom( V^{\otimes n }, V)\,, n \geq 1.$$

\vspace{3mm}

\noindent{\bf 1.3 Algebras over an operad.} Let $\oO$ be an operad in $\DGM$.  An {\it $\oO$-algebra}, $
(A,\rho)$,  is a differential graded module,
$A=\left(
\{A_i\}_{i\in \bZ}\,, d_A:A_i \to A_{i-1}\right)$, with an operadic
representation
$$
\rho :\oO \rightarrow \End _A \,.
$$
determined  by a sequence of maps differential graded
$\bk$-modules, called the {\it evaluation product}:
$$
\tilde \rho _n : \oO(n) \otimes A^{\otimes n} \to A \,,\quad \tilde \rho_n(
x\otimes a_1\otimes a_2\otimes ...\otimes a_n) =
\rho_n(x)( a_1\otimes a_2\otimes ...\otimes a_n).
$$
invariant under the action of $\Sigma _n$ and compatible with the composition
product of $\oO$.

$\oO$-algebras and homomorphisms of $\oO$-algebras is a category, denoted
$\oO\Alg$. If $\oO' $ is another operad in $\DGM $,
a homomorphism of operads
$ f :
\oO \to \oO'$ induces a natural functor $ f^\ast  : \oO'\Alg \to \oO\Alg$.

The {\it free $\oO$-algebra} generated by a differential graded module $V$ is
the differential graded module
$$ \displaystyle
F(\oO, V)= \bigoplus _{k=0}^\infty  \oO(k) \otimes _{\Sigma _k} V ^{\otimes
k}\,,
$$ with evaluation products  $\tilde \rho_n : \oO(n) \otimes F(\oO,V)^{\otimes
n} \to F(\oO,V) $   induced by the composition
products of $\oO$.   Any homomophism $ f : V \to W $
of graded modules  extends uniquely in homomorphism of graded modules
$F(\oO,f) : F(\oO,V) \to F(\oO,W) $. The  functor
$
F(\oO, -) : \DGM \to \oO\Alg
$
 is a left adjoint to the forgetful functor $ \oO\Alg \to \DGM$.

\vspace{3mm }

\noindent{\bf 1.4 The associative operad and the commutative operad.} For each $n \geq 0$, the canonical augmentation  $  \epsilon
_n  : \bk[\Sigma _n] \to \bk $   of the group
ring
$\bk[\Sigma _n]$ defines a homomorphism of operads in $\GM$
$$
\epsilon : \oA\to  \oC
$$
from the {\it associative operad}  $\oA$  to   the {\it commutative operad}
$\oC$ such that   $
\oA(n)=\bk[\Sigma _n]$ and $  \oC(n)=\bk$ with composition products given
respectively  by composition of permutations and
multiplication  in $\bk$.  An $\oA$-algebra is an associative differential
graded algebra while a $\oC$-algebra is a commutative differential   graded
algebra which is also, associative. Indeed, the
 representations $ \rho : \oA \to \End_A $ and $\oC \to \End _A$ are
the iterated products:
$$
\bk[\Sigma _n] \otimes_{\Sigma _n}
A^{\otimes n} =A^{\otimes n} \to A \mbox{ and } \bk \otimes_{\Sigma _n}
(A^{\otimes n}) =\left( A^{\otimes n}\right)_{\Sigma _n} \to A\,.  $$
 The  free
$\oC$-algebra (resp. a free $\oA$-algebra) generated  by the
differential graded module $V$ is the differential
graded module
$F(\oC,V)= \oplus  _{n=0}^\infty  \bk \otimes_{\Sigma _n}
V^{\otimes n}
=
\oplus_{n=0}^\infty   \left(  V^{\otimes n} \right)_{\Sigma _n} = S(V) $ with
graded commutative mutiplication of the
elements of $V$  (resp.$F(\oA,V)= \oplus  _{n=0}^\infty  \bk[\Sigma _n]
\otimes_{\Sigma _n}
V^{\otimes n}
=
\oplus_{n=0}^\infty     V^{\otimes n}  = T(V) $ the usual tensor algebra on
$V$).

\vspace{3mm }

\noindent{\bf 1.5} Let $A$ be an  $\oO$-algebra. Then $A$ is called
{\it almost free} if $A_\#= F(\oO_\#,V)$ for a graded module $V$. If the category of $\oO$-algebras is a closed model category (this is
the case whenever $\oO$ is a cofibrant operad \cite{[B-M]}) then any
$\oO$-algebra  admits a cofibrant model. An $\oO$-algebra is cofibrant if and only if it is a retract of an
almost free $\oO$-algebra.

\vspace{5mm }

{\bf 2. The  operads $\oE$ and $\oC_\infty$.}

\vspace{3mm }

\vspace{3mm }

\noindent{\bf 2.1}  Let us denote $\underline A =\left( \{A_n\}_{n\in \bN}\,,
d_i,
s_i\right)$ a simplicial  differential graded module with internal differential
$d_A : A_{n,q} \to A_{n, q-1}$.
Then the {\it total complex}  of the bicomplex
$$
\underline{A}_{p-1,q} \stackrel{ \sum (-1)^i d_i} {\leftarrow}
\underline{A}_{p,q} \stackrel{d_A}{\rightarrow}
\underline{A}_{p,q+1}
$$
is denoted by $\Tot_\ast ( \underline{A})$:
$$
\begin{array}{llll}
\Tot ( \underline{A}) =  \left( \{\Tot_n ( \underline{A}) \}_{n \in \bZ},
d\right)  &d :\Tot_n ( \underline{A})
 \to \{\Tot_{n-1} ( \underline{A})\\
\Tot_n ( \underline{A}) = \displaystyle \bigoplus _{p+q=n} \underline A _{p,q} &
dx = d_A x+(-1)^p \sum (-1)^i d_ix \,, x
\in \underline {A}_{p,q}\,.
\end{array}
$$

Let  $D_n (\underline A)$ be the subcomplex
generated by degeneracies in $ \underline A_n $. The
quotient complex $\Tot_\ast (\underline
A)/D_\ast(\underline A):=N_\ast(\underline A)$ is   the {\it normalized
differential graded module}. The quotient map
$\Tot (\underline A)
\twoheadrightarrow  N_\ast(\underline A )$ is a  chain equivalence,
\cite{[ML1]}-Theorem 6.1.

The singular chain complex of $X$ with coefficients in $\bk$  is
$C_\ast (X;\bk):= C_\ast(S_\ast(X;\bk))$  and the  normalized
chain complex is $N_\ast(X;\bk):=N_\ast (S_\ast(X)) $. The
subcomplex of normalized cochain complex $N^\ast(X,\bk)\subset
\Hom \left( C_\ast(X,\bZ),\bk \right)$ is, in this paper, simply
denoted $N^\ast(X)$. It follows from, \cite{[ML1]}-theorem 9.1
that $N^\ast$ is a contravariant functor from the category of
pointed topological spaces to the category of augmented
associative differential graded algebras.

The functor $N(-)$
$$\DGM^{\D^{op}}  \to \DGM  $$
is a  monoidal functor. The functor $N$ transforms   operads $\oO$
in $\DGM ^{\D ^{op}} $ (resp.  $\oO$-algebras in $\DGM ^{\D
^{op}}$ )  into  operads in $\DGM$  (resp.  agebras over an operad
in $\DGM$ ), \cite{[K-M]}-page 51. For further use we need the
slightly more general result.

\vspace{3mm} \noindent{\bf 2.2 Lemma.}{ \it Let $\underline A
=\{A_{p,q}\}_{p\in \bN,q\in \bZ}$ be a simplicial $\oO$-algebra.
Then the total complex $\Tot ( \underline{A})$ is an
$\oO$-algebra. Moreover, $N _\ast(\underline{A})$ is an
$\oO$-algebra and the quotient map $\Tot (\underline{A})\to  N
_\ast(\underline{A}) $ is an equivalence of $\oO$-algebras.}

\vspace{2mm}
\noindent{\bf Proof.}

 Let $\underline A$ and $\underline B$ be two simplicial graded modules  and
consider the {\it shuffle product}
$$
 \quad sh : C_p\underline A \otimes C_q\underline B \to  C_{p+q}(\underline
A\times \underline B)\,,  a\otimes b
\mapsto
\displaystyle  \sum _{\mu, \nu} (-1) ^{\epsilon (\mu)} s_\nu  a \times s_\mu b
\,,
\left\{
\begin{array}{rl}
a &\in A_p\\
 b&\in B_q
\end{array}
\right.
$$
when the sum is taken over the $p+q$ shuffles $ \mu _1 <\mu_2<...<\mu _p $,$ \nu
_1< \nu _2<...<\nu_q$, $\mu_i \,,\nu_j \in
\{1,2,...,p+q\}$,
$\epsilon (\mu ) $ is the graded signature of the  $(p,q)$-shuffle, $s_\mu =
s_{\mu _{p}}\circ
s_{\mu _{p-1}}\circ ...\circ s_{\mu _1}
$\,,
$s_\nu = s_{\nu _{q}}\circ s_{\nu _{q-1}}\circ ...\circ s_{\nu _1} $,
\cite{[ML1]}-Chapter 8. If we assume that
 $\underline A$ and $\underline B$ are  two simplicial differential graded
modules, one easily check that $sh $ commutes with
the differentials so that we obtain:
$$
 \quad sh : \Tot \underline A \otimes \Tot \underline B \to  \Tot (\underline
A\times \underline B)\,,
\mbox{ and }
\quad sh : N_p\underline A \otimes N_q\underline B \to  N_{p+q}(\underline
A\times \underline B)
\,.
$$
Since, when $A=B$, the
shuflle product is associative (and commutative) one defines the {\it iterated
shuffle product}
$$
sh^0=id\,, sh^{k+1} = (sh \otimes id)\circ sh ^k \,, \quad k \geq 0
$$

Let $\tilde \rho _{n,k} : \oO(k) \otimes (\underline A_n)^{\otimes k} \to A$ be
an evaluation product of $\underline A_n$. We
consider $\oO$ as a constant simplicial module and we define the map $ \hat \rho
_k : \oO(k) \otimes
\left( N_\ast (\underline A)\right) ^{\otimes k} \to N_\ast (\underline
A)
$  as the composite
$$
\begin{array}{lll}
\oO(k) \otimes  \left( N_\ast (\underline A)\right) ^{\otimes k}&=&  N_\ast
(\oO(k))\otimes N_\ast\left( \underline A)\right)
^{\otimes k}\\
&&\hspace{19mm} \downarrow  id \otimes sh ^k\\
&&  N_\ast (\oO(k))\otimes N_\ast\left( \underline A^{\otimes k}\right) \\
&& \hspace{19mm} \downarrow  sh \\
 N_\ast (\underline A) &\stackrel{ N_\ast (\tilde \rho _k)}{\leftarrow} &
 N_\ast \left(\oO(k))\otimes \underline A^{\otimes k}\right)
\end{array}
\,.
$$

It is then easy to check that the $\hat \rho _n$ are composition products.

\hfill $\square$

\vspace{2mm} \noindent{\bf 2.3} {\bf The operad $\oE$.} The operad
$\oE$ (also called the Barrat-Eccles operad, \cite{[B-F]}-1.1.) is
an operad in the category of differential graded $\bk$-modules
such that:
$$
\oE(n)= N_\ast \left( \W (\Sigma _n)\right) = \mbox{ \rm the normalized bar
construction on } \Sigma _n .
$$

\vspace{2mm}
\noindent{\bf 2.4}  {\bf  Some properties of $\oE$.}

1) The operad $\oE$ is a resolution of the operad $\oC$: the homomorphism of
operads
 $
\overline \epsilon : \oE \stackrel{\sim}{\twoheadrightarrow} \oC
$ is
defined by the augmentations of the bar resolution for each component.

2)  $\oE$ is an $E_\infty$-operad. Recall that an operad $\oO=\{\left(\oO(n)\right)_i\}_{i\geq 0}$ in $\DGM$  is an
{\it
$E_\infty$-operad } if each $\oO(n)$ is an acyclic $\Sigma _n $-free module.

3)  The natural map
$\epsilon :
\oA
\to
\oC$ factorises as
$
\oA \to \oE \stackrel{\overline{\epsilon}}{\longrightarrow} \oC
$. In particular, $\oE$-algebras are associative algebras.

4) The operad $\oE$ is not cofibrant.

5) Berger and Fresse \cite{[B-F]} have proved that the normalized
singular cochains
$N^*(-)$ is a functor from the category of
topological spaces to the category of $\oE$-differential algebras.

6)  If $\bk =\bF _p$ is the prime field
of characteristic
$p$   and $A$ is a $\oE$-algebras   then $H(A)$  is   an unstable algebra over
the big Steenrod
algebra (\cite {[Ma]}, \cite{[C-L]}). Indeed, consider the Standard small free resolution ${\cal W}$ (resp.
${\cal W}'$ of the cyclic group of order $p$;
$\pi
\subset
\Sigma _p$ (resp. of $\Sigma _p$. Thus we obtain the homomorphism  $ {\cal W}
\to {\cal W}' \to \oE (p) $ and  the
evaluation product
$\tilde \rho _p :
\oE (p)
\otimes A^p
\to A$ restricts to the structure map $ {\cal W} \otimes A^{\otimes p} \to A$
considered by May, \cite{[Ma]}, in order to
define ``algebraic Steenrod operations''. Recall that the big Steenrod algebra,
denoted $B_p$, is such that the quotient
$B_p/(P^0=id)$ is the usual Steenrod algebra, see \cite{[Ma1]}-theorem 1.4. In
particular, Adem relations are satisfied (\cite{[Ma]} and \cite{[C-L]}).

\vspace{3mm}
\noindent{\bf 2.5} {\bf The operad $\oC_\infty$.}  Let  $\oC_\infty$ be a cofibrant replacement of $\oC$. There exists
a quasi-isomorphism of operads
$$
\oC_\infty \to \oE .
$$

In particular,  by  remarks 5) and 6) above, $N^\ast (X)$ is a
$\oC_\infty $-algebra  and any quasi-isomorphism of $\oC_\infty
$-algebras $ A \to N^\ast(X)$ identifies the action of the large
Steenrod operations on $H(A)$ to the usual action on $N^\ast (X)$.

By (1.5)  any $\oC_\infty $-algebra admits a cofibrant replacement which is
an almost free $\oC_\infty$-algebra, see \cite{[H-JPAA]} and \cite{[B-M]}

\vspace{5mm }

{\bf 3. The Operadic Hochschild chain complex.}

\vspace{3mm }
\noindent{\bf 3.1} Let us recall that  the category of $\oO$-algebras has all
limits and colimits,
\cite{[G-J]}-Theorem 1.13. In particular, \cite{[Ma1]}-3,  the coproduct of two
$\oO$-algebras $A$ and $B$ is an
$\oO$-algebra, denoted
$A\coprod B$. Hereafter, we will use the following notation:
$
A \stackrel{l}{\to} A\coprod B \stackrel{r}{\leftarrow }  B
$
for the natural inclusions and  $ A\coprod A \stackrel{\nabla}{\to} A$ for the
folding map.  The
symmetric group
$\Sigma _n$ acts on $A^{^\coprod n}= A\coprod ...\coprod A $ ($n$ terms) by
permutations of factors. We denote by $\tau_n :
A^{^\coprod n} \to A^{^\coprod n} $  the homomophism corresponding to  the
cyclic permutation $(n,1,2,...,n-1) \in \Sigma
_n$.

\vspace{3mm }
\noindent{\bf 3.2}  Let $\oO$ be an operad in $\DGM$ and $A$ be an $\oO$-algebra.
We denote by $\underline A$ the simplicial
$\oO$-algebra
$$
\begin{array}{lllll}

\underline{A}_n =  A^{ ^\coprod n+1}\,, n \geq 0\,,\quad  d_i : A^{^\coprod n+1}
\to A^{^\coprod n} \,,\quad  s_i :
A^{^\coprod n}
\to A^{^\coprod n+1}\\
d_i = \left\{
\begin{array}{ll}
id ^{^\coprod i}\coprod \nabla \coprod id ^{^\coprod n-i}& \mbox{ if } i=
0,1,...,n-1\\
\left( \nabla \coprod id ^{^\coprod n}\right) \circ \tau _n &\mbox{ if } i= n
\end{array}
\right. \,, \quad
s_i = id ^{^\coprod i}\coprod l \coprod id ^{^\coprod n-i}
\,.
\end{array}
$$

The normalization $N(\underline A)$, (see 2.2), of the simplicial $\oO$-algebra
$\underline{A}$ is an $\oO$-algebra, denoted
$\hC (\oO, A)$, and called the {\it Hochschild chain complex of the
$\oO$-algebra $A$}. The homology of the $\hC (\oO, A)$ is
the {\it operadic Hochschild homology}, denoted $\hH (\oO, A)$.

\vspace{3mm } \noindent{\bf 3.3}  If
$(A,d_A)$ is any {\it associative} differential graded algebra
supposed  unital and augmented we denote by
$\overline A$ the kernel of the augmentation. The (classical)
Hochschild chain complex is defined as follows:
$$
 {\frak C}_*A = \{ {\frak C}_kA\}_{k \geq 0} \,, \quad  {\frak C}_kA   = A
\otimes s {\overline A }^{\otimes k}\,,
$$
with $A= \bk\oplus \overline A$. A  generator of  ${\frak C}_kA$ is written
$a_0[sa_1|sa_2|...|sa_k]$ if $k>0$
and  $a[\,]$ if $k=0$. We set $\epsilon _i = |a_0|+
|a_1|+|sa_2|+...+|sa_{i-1}|$, $i\geq 1$. The differential $d=d^1+d^2$
is defined by:
$$
\begin{array}{ll}
d^1 a_0[a_1|a_2|...|a_k] &= da_0[a_1|a_2|...|a_k] - \sum _{i=1}^k (-1)^{\epsilon
_i} a_0[a_1|...|da_i|...|a_k]\\

d^2 a_0[a_1|a_2|...|a_k] &=  (-1) ^{|a_0|} a_0a_1[a_2|...|a_k] +
 \sum _{i=2}^k (-1)^{\epsilon _i} a_0[a_1|...|a_{i-1}a_i|...|a_k]\\
& \qquad -  (-1)^{|sa_k|\epsilon _k} a_ka_0[a_1|...|a_{k-1}]
\end{array}
$$
 Consider the shuffle map, \cite{[Lo]} 4.2.1,
       $sh: {\frak C}_* A\otimes {\frak C}_* A
       \longrightarrow {\frak C _*
       }(A\otimes
       A)$ defined by:
       $$sh(a_0[ a_1 \vert a_2\vert ... \vert a_n] , b_0[ b_1
       \vert b_2\vert ... \vert
       b_m] ) = (-1)^t \sum_{\sigma\in \Sigma_{n,m}}
       (-1)^{\epsilon (\sigma)}a_0\otimes b_0[
       c_{\sigma(1)}\vert ...\vert c_{\sigma (m+n)}]
       $$
       where
       $t = \vert b_0\vert (\vert a_0\vert +...+\vert a_{n}\vert)$,
       $
       \quad \quad
       c_{\sigma (i)} = \left \{ \begin{array}{ll}
       a_{\sigma (i)}\otimes 1, &1\leq i \leq n \\
       1\otimes b_{\sigma (i-n)}, &n+1\leq i \leq n+m
       \end{array}
       \right.
       $\\
       and $\epsilon (\sigma) = \sum \vert c_{\sigma (i)}\vert \vert c_{\sigma
       (m+j)}\vert$, summed over all pairs
       $(i, m+j)$ such that $\sigma (m+j) <\sigma (i)$.
       Clearly, $sh$ induces
a chain map still denoted $sh$, ${\frak C}_*
       A\otimes {\frak C}_* A \longrightarrow {\frak C
       _* }(A\otimes
       A)$.
       Let $HH_\ast(A)$ be the homology of ${\frak C}_* A$.

 If $A$ is commutative (in the graded sense) then the multiplication
       $\mu_A :A\otimes A \rightarrow A$
       is a homomorphism of differential graded algebras. Thus the composite
       ${\frak C _* }\mu_A \circ sh: {\frak C _* }A\otimes
       {\frak C _* }A \rightarrow {\frak C _* }A$ defines a multiplication on
       ${\frak C _*}A$ which
       makes it into a commutative differential graded  algebra \cite{[Lo]}
4.2.2.

If $A$ is an associative (non commutative) differential graded algebra, there is
no interesting product
on
${\frak C}_\ast (A)$ while $\hC(\oA,A)$ is naturally an associative differential
graded algebra. Obviously,
${\frak C}_\ast (A)_\# \neq \hC(\oA,A)_\#$.  N. Bitjong and second named author,
\cite{[BT1]},  have  proved that  for any
strongly homotopy  commutative $\bk$-algebra $A$, in the sense of
\cite{[Mu]}, there is a well defined   cup product on
$HH_\ast(A)$ , (\cite{[BT1]}-Theorem 1), which is induced from a non canonical
product on ${\frak C}_\ast (A)$.
  In   the formalism of operads, a strongly homotopy  commutative $\bk$-algebra
$A$ is an associative $\oB_\infty $-algebra, in the sense of
\cite{[G-J]}- 5.2, (see \cite{[BT1]}-Proposition 2). The graded
vector space $\hC(\oB_\infty, A)$ is not isomorphic to ${\frak
C}_\ast (A)$. An interesting question is:
 {\it  Let $A$ be a strongly homotopy  commutative
algebra $A$  does  $\hH (\oB_\infty ,A) \cong  HH_\ast(A)$ as commutative graded
algebras?}

\vspace{3mm }
\noindent{\bf 3.4} Let $A$ be an associative (unital) differential
graded algebra.  There is classically associated to $A$ an other
simplicial differential  graded algebra, which we denote
$\underline{\underline{A}}$ and which is defined as
follows (see
\cite{[Jo]}-Exemple 1.4):
$$
\begin{array}{lllll}
\underline{\underline{A}}_n =  A^{ \otimes  n+1}\,, n \geq 0\,,\quad  d_i :
A^{\otimes n+1} \to A^{\otimes n} \,,\quad  s_i
: A^{\otimes n}
\to A^{\otimes n+1}\\
d_i = \left\{
\begin{array}{ll}
id ^{\otimes i}\otimes \mu_A \otimes id ^{\otimes n-i}& \mbox{ if } i=
0,1,...,n-1\\
\left( \otimes \mu_A \otimes id ^{\otimes  n}\right) \circ \tau _n &\mbox{ if }
i= n
\end{array}
\right. \,, \quad
s_i = id ^{\otimes i}\otimes 1 \otimes id ^{\otimes n-i}
\,.
\end{array}
$$
where $\mu _A$ demotes the multiplication of $A$ and $\tau _n $ the map $
a_0\otimes a_1\otimes...\otimes a_n \mapsto
(1)^{(|a_0|+...+|a_{n_1}|)|a_n]} \to  a_n\otimes a_0\otimes...\otimes a_{n-1}$.

The complex $\tilde  {\frak C}(A)= N_\ast
\underline{\underline{A}}$ is the {\it unreduced Hochschild chain
complex}. By \cite{[ML1]}-Chapter X-Corollary 2.2 and Theorem 9.1, the
maps
$  (id \otimes s^{\otimes n}) $
define a quasi-isomorphism
$$
 (1)\qquad  \tilde  {\frak C}(A)  \to {\frak C}_* (A)\,.
$$

If we assume that $A$ is commutative, $\tilde  {\frak C}_\ast (A)$ is also a
differential graded algebra.  The multiplication
is the   shuffle product defined in the same  way that the  shuffle product on
${\frak C}_* (A) $. Therefore the quasi-isomorphism (1) is a homomorphism of
differential graded algebras.

\vspace{3mm }
\noindent{\bf 3.5}  { \bf Theorem A.} {\it Assume that  A is associative
commutative  differential graded algebra.   There
exists a natural  quasi-isomorphism  of commutative differential graded algebras
$$
\hC(\oC,A) \to {\frak C}(A)\,.
$$}

{\bf Remark.} If the $\oC$-algebra $A$ is considered as a
$\oC_\infty$-algebra then, by naturality there is a surjective
homomorphism of $\oC_\infty$-algebras
$$
\hC(\oC_\infty,A)\to\hC(\oC,A)\,.
$$
Does this map induces an isomorphism in homology?

\vspace{2mm} {\bf Proof.}  Let $B$ be a  $\oC$-algebra. By universal property, there exists a
natural isomorphism
$\Phi _{A,B} : A\coprod B \to A \otimes B$  of commutative differential graded
algebras such that the following diagrams
commute where we put  $\Phi _{A,A}=\Phi_1$.
$$
\begin{array}{lll}
A\coprod A &\stackrel{\nabla} {\to} & A\\
\Phi _1 \downarrow && || \\
A\otimes A &\stackrel{\mu} {\to} & A
\end{array}
\,,
\quad
\begin{array}{lll}
A\coprod A &\stackrel{T} {\to} & A\coprod A\\
\Phi _1 \downarrow && || \\
A\otimes  A &\stackrel{T} {\to} & A\otimes  A
\end{array}
$$
where, $\mu $ denotes the usual product on $A\otimes A$
and $T$ the usual twisting map. The
associativity properties permit
iteration so that we obtain for any $n\geq 0$
an isomorphism $ \Phi _ n : A ^{^\coprod n_1}\to A^{\otimes n+1}$. These
$\Phi_n$'s induce an isomorphism  of simplicial
differential graded modules $
\underline A \to \underline{\underline {A}}$ which in turn induces an
isomorphism
$$
(2) \qquad  \hC(A)= N_\ast  (\underline A) \to N_\ast (
\underline{\underline{A}})=\tilde{\frak C}_\ast  (A)\,.
$$
Composition of the isomorphism (2) with the quasi-isomorphism (1) gives the
quasi-isomorphism
$$
(3) \qquad  \hC(A) \to {\frak C}_\ast  (A)  \,.
$$
 On the other hand,   the $\oC$-algebra structure on $ \hC(A)$ is such that the
isomorphism (2) is an
isomorphism of differential graded algebras. Thus (3) is a quasi-isomorphism of
commutative differential graded algebras.

\hfill{$\square$}

\vspace{3mm } \noindent{\bf 3.6}  An operad $\oC_\infty $ is a
Hopf operad "up to homotopy" (for the notion of Hopf operad we
refer to \cite{[G-J]}-5.3), that is to say it has a diagonal wich
is not coassociative but only up to homotopy). In this case, the
tensor product of two $\oC_\infty $-algebras $A$
 and $B$  is a $\oC_\infty $-algebra \cite{[C]} with
underlying  differential graded module being the tensor product of
the underlying differential graded modules, denoted $A \otimes B$.
Indeed, Hinich \cite{[H-JPAA]}  has  proved  that if the
$\oC_\infty$-algebras, $A$ and $B$ are cofibrant there exists a
natural  quasi-isomorphism
$$
\Phi_{A,B} : A\coprod B \to A\otimes B\,.
$$

Let $ f : A \to A' $ and $ g: B \to B' $ be two homomorphisms of $\oC_\infty
$-algebras. Then, by naturallity of $\Phi
_{A,B}$,  we obtain the commutative diagram
$$
\begin{array}{lll}
A \coprod B &\stackrel{ f\coprod g}{\to} & A' \coprod B'\\
\Phi _{A,B}  \downarrow && \downarrow \Phi _{A',B'}\\
A \otimes  B &\stackrel{ f \otimes  g }{\to} & A' \otimes  B'
\end{array}
$$
If we assume that $A$ and $A'$ are cofibrant and that $f$ and $g$ are
quasi-isomorphisms, then by \cite{[Ma1]}-Theorem 3.2,
the homomorphism $ f\coprod g$ is a quasi-isomorphism. Therefore, $f\otimes g$ is
also a quasi-isomorphism.

\vspace{3mm }
\noindent{\bf 3.7} Assume that $ A= F(\oO, V)$  is an almost free $\oO$-algebra.
Thus, we have the sequence of direct summands
$$
\underline A_n \supset \oO(k)\otimes_{\Sigma _k} (V^{\oplus n+1})^{\otimes k }\supset \oO(k)
\otimes_{\Sigma _k}\left( V^{\otimes k_0}\otimes V^{\otimes
k_1}\otimes...\otimes  V^{\otimes k_n}\right)\,,
$$
with $k=k_0+k_1+...+k_n$. Therefore, an element of $\underline A_n
$ is finite sum of elements of the form $ x \otimes
v_1....v_{k_0}\otimes v_{k_0+1}...v_{k_0+k_1}\otimes ... \otimes
v_{k_0+...+k_{n-1}+1}...v_k$ with $x \in \oO(k)$ and $v_{k_0+
...+k_{i-1}+1}...v_{k_0+ ...+k_{i}} \in V ^{\otimes k_i}$  with
usual  convention $V^{\otimes 0}=\bk$. With this notation we
obtain explicit formulas for the map $r$, $l$ and $\nabla $
defined in 3.1:
$$
\begin{array}{lll}
r : \underline A_0 \to \underline A_1 & r( x \otimes v_1....v_{k_0} = x \otimes
v_1....v_{k_0}\otimes 1\\
l : \underline A_0 \to \underline A_1 & l( x \otimes v_1....v_{k_0} = x \otimes
1\otimes v_1....v_{k_0}\\
\nabla  : \underline A_1 \to \underline A_0 &
\nabla (x \otimes v_1....v_{k_0}\otimes v_{k_0+1}...v_{k_0+k_1})= x \otimes
v_1....v_{k_0} v_{k_0+1}...v_{k_0+k_1}\,.
\end{array}
$$
Therefore,
$$
\begin{array}{l}
d_i : \underline A_{n+1}  \to \underline A_n  \,, \quad
d_i\left( x \otimes v_1....v_{k_0}\otimes v_{k_0+1}...v_{k_0+k_1}\otimes ...
\otimes v_{k_0+...+k_{n-1}+1}...v_k\right)= \\
\qquad \left\{
\begin{array}{ll}
x \otimes v_1....v_{k_0}\otimes ... \otimes  v_{k_0 + k_i+1}...v_{k_0+...+
k_{i+2}}\otimes ... \otimes
v_{k_0+...+k_{n-1}+1}...v_k &\mbox{ if } i=0,1,...,n-1\\
=
x \otimes  v_{k_0+...+k_{n}+1}...v_k  v_1....v_{k_0 }\otimes ... \otimes
v_{k_0+...+k_{n-2}+1}...v_{k_0+...+k_{n-1}} &\mbox{ if } i=n\\
\end{array}
\right.\\
s_i: \underline A_{n}  \to \underline A_{n+1} \,, \quad
s_i \left(x \otimes v_1....v_{k_0}\otimes v_{k_0+1}...v_{k_0+k_1}\otimes ...
\otimes v_{k_0+...+k_{n-1}+1}...v_k\right)=\\
\qquad x \otimes v_1....v_{k_0}\otimes ...\otimes 1 \otimes ...\otimes
v_{k_0+...+k_i+1 }...v_{k_0+...+k_{i+1} }\otimes   ...
\otimes v_{k_0+...+k_{n-1}+1}...v_k\,.
\end{array}
$$

Let us denote by $\underline A ^+ _n$ the submodule of the $\bk$-module
$\underline{A}_n$ generated by the elements of the
form $ x \otimes v_1....v_{k_0}\otimes v_{k_0+1}...v_{k_0+k_1}\otimes ...
\otimes v_{k_0+...+k_{n-1}+1}...v_k$ with $x \in
\oO(k)$ and $v_{k_0+ ...+k_{i-1}+1}...v_{k_0+ ...+k_{i}} \in V ^{\otimes k_i}$
such that each $k_i>0$. The above fomulas for
$d_i$ and $s_i$ show that

a) the  graded module $\underline A^+_n$ is stable for the $d_i$'s but not
stable for the $s_i$'s.

b) the submodule $D\underline {A} _n$ of $\underline {A} _n$
generated by all degenerate elements ($D\underline {A} _0= 0$) is
exactly the submodule generated by the  elements $x \otimes
v_1....v_{k_0}\otimes v_{k_0+1}...v_{k_0+k_1}\otimes ... \otimes
v_{k_0+...+k_{n-1}+1}...v_k$ such that at least one $k_i=0$.

Since $ \hC(\oO, A)= N(\underline A)= \Tot (\underline A)/D(\underline A)$ and
since $d_A(\underline{A}^+_n) \subset
\underline{A}^+_n$  we have proved the first part of the next result.

\vspace{3mm } { \bf Theorem B.} {\it Assume that $ A= F(\oO, V)$
is an almost free $\oO$-algebra. Then  the restriction of the
natural chain equivalence $\Tot \underline A \to N_\ast \underline
A$  to the total complex, $\Tot ( \underline{A}^+)$  is an
isomorphism of  differential graded modules
$$
\Tot ( \underline A^+) \to \hC(\oO,A)\,.
$$ Moreover, this homomorphism is an
isomorphism of $\oO$-algebras.}

\vspace{2mm} {\bf End of proof.} Let us precise first that the evaluation
products of the simplicial algebra $\underline A$
are determined by the maps
$$
\oO(k)_q \otimes  \underline{A}^+_{p_1,q_1}\otimes
\underline{A}^+_{p_2,q_2}\otimes ... \otimes
\underline{A}^+_{p_k,q_k} \to
\underline{A}^+_{p_1+p_2+...+p_k,q+q_1+...q_1}
$$
which are  explicitely given by the shuffle products and the evaluation product
of $\oO$ (see proof of lemma 2.2). This
implies that  the following diagram commutes
$$
\begin{array}{lll}
\oO(n)_q \otimes \underline{ A} ^+_{p_1,q_1} \otimes ....\otimes \underline{ A}
^+_{p_n,q_n } & \hookrightarrow &
\oO(n)_q \otimes \underline{ A} _{p_1,q_1} \otimes ....\otimes \underline{ A}
_{p_n,q_n }\\
\hspace{12mm} \downarrow & & \hspace{12mm}\downarrow \\
  \underline{ A}^+ _{p_1+...+p_n,q+q_1+...+q_n}
 & \hookrightarrow &
\underline{ A} _{p_1+...+p_n,q+q_1+...+q_n}
\end{array}
$$
Therefore, each $\underline A^+_n$ is a sub $\oO$-algebra of $\underline A_n$.

\hfill{$\square$}

It results from the formula above that $d_0=d_1 : \underline{A}_1 \to
\underline{A}_0$ and that $ \displaystyle d_i \left(
\bigoplus _{p>0, q
\in
\bZ}
\underline{A}_{p,q}\right) \subset \bigoplus _{p>0, q
\in
\bZ}
\underline{A}_{p,q}$. Thus we obtain:

\vspace{3mm }
{ \bf Proposition.} {\it Assume that $ A= F(\oO, V)$  is an almost free
$\oO$-algebra. Then we have the
natural splitting of $\oO$-algebras}
$$
    \Tot ( \underline{A}^+) = (A, 0)  \oplus \left( \bigoplus_{p>0}
\underline{A}^+ _{p,q}, d\right)  \,.
$$

\vspace{5mm }

{\bf 4. Free loop space.}

\vspace{3 mm}
\noindent {\bf 4.1}   Write  ${\bf Top}$ (resp. ${\bf Costop}$) for  the
category of topological spaces (resp. of cosimplicial topological spaces).
The geometric realization of a cosimplicial set  is the  covariant functor
$$
\vert \vert . \vert \vert: {\bf
Costop}\rightarrow {\bf Top}
 \,, \quad \underline Z \mapsto  \vert \vert \underline {Z}\vert \vert =
{\bf Costop}(\underline{\triangle}, \underline{Z})
\subset \prod_{n\geq 0} {\bf Top}(\triangle^n,{\underline Z}(n))\,,
$$
where $
 \vert \vert \underline {Z}\vert \vert$ is
equipped with the
topology induced by this inclusion. Here $\underline{\triangle}$ denotes
the cosimplicial space defined by
$\underline{\triangle}(n)
= \triangle^n$
 with the usual coface and codegeneracy maps. If $\underline Z$ is any  cosimplicial topological space,
 then $N^\ast\underline {Z}$ is a simplicial cochain complex with total complex
$\Tot (N^\ast\underline {Z})$:
$$
(\Tot_n ( N^*\underline {Z}) =
\oplus_{p-q=n} N^q\underline{Z}(p) \,, \quad Dx =
\sum^p_{i=1}(-1)^iC^\ast(d_i)+(-1)^p\delta x
\,, x\in
N^\ast \underline Z(p)
\,.$$
 The
 $d_i$ are the coface operators of $\underline Z$  and $\delta$ is the internal differential of
$N^\ast (\underline Z(p))$. Recall that in general the natural
map $ \Tot (N^\ast\underline {Z}) \to N^\ast(
||\underline Z||)$ is not a weak equivalence, \cite{[BS]}. It results
from  lemma 2.2 and 2.3-5 that the total  complex $\Tot
(N^\ast(\underline {Z}))$ is naturally a
$\oC_\infty$-algebra and $ \Tot (N^\ast(\underline {Z})) \to N^\ast(
||\underline Z||)$ a homomorphism of $\oC_\infty$-algebras.

Hereafter we denote $N^\ast\left(N^\ast (\underline Z)\right)$ the normalization
of $\Tot \left(N^\ast (\underline
Z)\right)$.

\vspace{3 mm}
\noindent {\bf 4.2}  One of the interest for considering cosimplicial spaces is the following  result,  \cite{[BS]}-Proposition 5.1, (see also
\cite{[P-T]}-Corollary 1): {\it If $\underline L $ is a simplicial set and $T$ a topological space then
the cosimplicial space  $T ^{\underline L}$  is such that  there is a
homeomorphism:
$$ \vert \vert T^{\underline L}\vert \vert :=
 {\bf Costop}(\underline{\triangle},  T^{\underline L}) \cong  {\bf Top}(\vert\underline L \vert , T) = T^{\vert\underline L\vert} \,.$$
}
In particular, if  we consider the simplicial set
$K$ defined as follows:
$K(n) = {\Bbb Z}/(n+1){\Bbb Z}, $ and, if $\overline {k}^{n}$ denotes an
element in ${\Bbb Z}/n{\Bbb Z}$,
the face maps  $d_i: K(n)\rightarrow K(n-1)$ with $0\leq i\leq n-1$ and the
degeneracy maps  $s_j:  K(n)\rightarrow K(n+1)$
with
$0\leq j\leq n$
are:
 $$
d_i \overline { k}^{n+1}=\left\{
\begin{array}{rl}
\overline {k}^{n} & \mbox{ if } k\leq i\\
\overline {k-1}^{n} & \mbox{ if } k>i
\end{array}
\right.
\qquad
 s_j\overline {k}^{n+1}=
       \left \{
       \begin{array}{rl}
        \overline { k}^{n+2}  &\mbox{ if } k\leq j\\
       \overline {k+1}^{n+2} &\mbox{ if }  k>j .
       \end{array}
       \right.
$$
 and $ d_n \overline {k}^{n+1}  = \overline {k}^{n}$. The  geometric realization
of $K$, \cite{[BF]} (proposition 1.4),
$\vert K\vert$ is homeomorphic to the circle $S^1$. Therefore,  the cosimplicial model,  $\underline X$,   for the free loop space,
used by
Jones, \cite{[Jo]},
$$
 \begin{array}{ll}
\underline{X}(n) & = Map (K(n), X) = \underbrace {X\times....\times X}_{(n+1)
\mbox{-times}}\\
d_i(x_0,x_1,...,x_n) &= (x_0,x_1,...,x_i,x_i, ... x_n) \,, 0\leq i \leq n\\
d_{n+1}(x_0,x_1,...,x_n) &= (x_0,x_1,...,x_n, x_0)\\
s_j(x_0,x_1,...,x_n)&= (x_0,x_1,...x_j,x_{j+2},...,x_n)\,, 0\leq j \leq n.
\end{array}
$$
is such that  $ \vert \vert \underline X \vert \vert \cong  {\bf Top}(\vert
\underline K
\vert , X) = X^{S^1}$, ($\cong $ means homeomorphism). From lemma 2.2,
we deduce then:

\vspace{3mm} \noindent{\bf Proposition.} {\it If $X$ is simply
connected, the natural map $\Tot (N ^ \ast {\underline X}) \to
N^\ast(X^{S 1})$ is a quasi-isomomorphism of
$\oC_\infty$-algebras.}

\vspace{3 mm}
\noindent {\bf 4.3} {\bf Theorem C.} {\it Let $X$ be a
simply connected space such that each $H^i(X)$
is finitely generated. Given an almost free model of the space $X$
$$
\varphi _X : M_X= \left( F(\oC, V), d\right) \to N^\ast (X)
$$
 there exists a   natural quasi-isomorphism  of
$\oC_\infty$-algebras
$$
\hC(\oC_\infty, M_X) \longrightarrow  N^\ast\left(N^\ast (\underline X)\right) \,.
$$}

\vspace{2mm}
\noindent{\bf Proof.}  Let $\varphi _X : M_X= \left( F(\oC, V), d\right) \to
N^\ast (X)$  (resp.
 $\varphi _Y : M_Y= \left( F(\oC,
V), d\right) \to N^\ast (Y)$) be a almost free model for the space $X$ (resp.
for the space $Y$). By universal property, we obtain the commutative diagram
$$
\begin{array}{cccccccccllll}
M_X &\stackrel{l}{\rightarrow} M_X &\coprod M_Y& \stackrel{r}{\leftarrow} &M_Y\\
\sim \downarrow \varphi_X && \hspace{-10mm} \psi_{X,Y} \downarrow && \sim
\downarrow \varphi_Y \\
N^\ast(X)& \stackrel{N^\ast(pr_X)}{\rightarrow}& N^\ast(X\times Y)
&\stackrel{N^\ast(pr_Y)}{\leftarrow}
&N^\ast(Y)\\
\end{array}
$$
where $pr_X$ and $pr_Y$ (resp. $i_X$ and $i_Y$) are the natural projections
(resp. inclusions). We have also the following
commutative diagrams
$$
\begin{array}{cccccccccllll}
M_X \coprod M_X &\stackrel{\nabla}{\rightarrow}& M_X \\
 \psi_{X,X} \downarrow && \sim \downarrow \varphi_X \\
 N^\ast(X\times X) &\stackrel{N^\ast(\Delta_X)}{\rightarrow} &N^\ast(X)\\
  \quad
\end{array}
\begin{array}{cccccccccllll}
M_X \coprod M_X &\stackrel{T^{^\coprod}}{\rightarrow}& M_X \coprod M_X \\
 \psi_{X,X} \downarrow && \downarrow \psi_{X,X}\\
 N^\ast(X\times Y) &\stackrel{N^\ast(T_X)}{\rightarrow} &N^\ast(X)\\
\end{array}
$$
where $\Delta _X $ is the diagonal  and $T_X$ the topological
interchange map.

By \cite{[Ma1]}-Lemma 5.2, we know that if each $H^i(X)$
is finitely generated then  $\psi_{X,X} $ is a weak equivalence of differential graded
module. It is now easy to prove that
iteration furnishes a  quasi-isomorphism of simplicial
$\oC_\infty$-algebras $
\underline M_X \to N ^\ast (\underline X)
$.

\hfill{$\square$}

{\bf Acknowledgements.} The authors would like to thank Benoit
Fres\-se for pointing out a mistake in a preliminary version and
Martin Markl for his careful reading and suggestions.

\vspace{1cm}

\hspace{-1cm}\begin{minipage}{19cm}
\small
\begin{tabular}{lll}
dchataur@crm.es
&$\,$\hspace{4cm}$\,$&jean-claude.thomas@univ-angers.fr\\
CRM Barcelona&&
D\'epartement de math\'ematique \\

Institut d'Estudis Catalans &&  Facult\'e des
Sciences \\

Apartat 50E  && 2, Boulevard
Lavoisier     \\

 08193 Bellaterra, Espagne   &&
 49045
Angers, France

\end{tabular}

\end{minipage}

\end{document}